\documentclass{amsart}
\usepackage{amssymb, latexsym}
\usepackage{enumerate}

\newcommand \maxJ {\max{(1,|J|)}}

\newcommand \ulJ{u^{l,J}}
\newcommand \unlJ{u^{nl,J}}

\theoremstyle{plain}
\newtheorem{thm}{Theorem}

\newtheorem{rem}[thm]{Remark}
\newtheorem{prop}[thm]{Proposition}

\numberwithin{equation}{section} \numberwithin{thm}{section}

\begin{document}

\title[Defocusing cubic wave equation]{ Adapted Linear-Nonlinear Decomposition and Global well-posedness for solutions \\
to the defocusing cubic wave equation on $\mathbb{R}^{3}$ }
\author{Tristan Roy}
%\thanks{}
\address{University of California, Los Angeles}
\email{triroy@math.ucla.edu}
%\subjclass{35Q55} \keywords{nonlinear Schr\"odinger equation,
%well-posedness}

\vspace{-0.3in}

\begin{abstract}
We prove global well-posedness for the defocusing cubic wave
equation

\begin{equation}
\left\{
\begin{array}{ccl}
\partial_{tt} u  - \Delta u & = & -u^{3} \\
u(0,x)& = & u_{0}(x)    \\
\partial_{t} u(0,x) & = & u_{1}(x)
\end{array}
\right. \nonumber
\end{equation}
with data $\left( u_{0}, \, u_{1} \right) \in H^{s} \times H^{s-1}$,
$1 > s > \frac{13}{18} \simeq 0.722$. The main task is to estimate
the variation of an almost conserved quantity on an arbitrary long
time interval. We divide it into subintervals. On each of these
subintervals we write the solution as the sum of its linear part
adapted to the subinterval and its corresponding nonlinear part.
Some terms resulting from this decomposition have a controlled
global variation and other terms have a slow local variation.
\end{abstract}

\maketitle

\section{Introduction}

We shall study the defocusing cubic wave equation on
$\mathbb{R}^{3}$

\begin{equation}
\begin{array}{ll}
\partial_{tt} u  - \Delta u & =  -u^{3} \\
\end{array}
\label{Eqn:WaveEq}
\end{equation}
with data $u(0)=u_{0}$, $\partial_{t}u(0)=u_{1}$  lying in $H^{s}$,
$H^{s-1}$ respectively. Here $H^{s}$ is the standard inhomogeneous
Sobolev space i.e $H^{s}$ is the completion of the Schwartz space
$\mathcal{S}(\mathbb{R}^{3})$ with respect to the norm

\begin{equation}
\begin{array}{ll}
\| f \|_{H^{s}} & :=  \| (1+ D)^{s} f \|_{L^{2}(\mathbb{R}^{3})}
\end{array}
\end{equation}
where $D$ is the operator defined by

\begin{equation}
\begin{array}{ll}
\widehat{Df}(\xi) & := |\xi| \hat{f}(\xi)
\end{array}
\end{equation}
and $\hat{f}$ denotes the Fourier transform
\begin{equation}
\begin{array}{ll}
\hat{f}(\xi) & := \int_{\mathbb{R}^{3}} f(x) e^{-i x \cdot \xi} \,
dx
\end{array}
\end{equation}
We shall focus on the strong solutions of the defocusing cubic wave
equation on some interval $[0,T]$ i.e real-valued maps $u$,
$\partial_{t} u$ that lie in  $C \left([0, \, T], \,
H^{s}(\mathbb{R}^{3}) \right)$, $C \left( [0, \,T], \, H^{s-1} (
\mathbb{R}^{3}) \right)$ respectively and that satisfy for $t \in
[0, \, T] $ the following integral equation

\begin{equation}
\begin{array}{ll}
u(t) & = \cos(tD) u_{0} + D^{-1} \sin(tD) u_{1} - \int_{0}^{t}
D^{-1} \sin \left( (t-t^{'}) D \right) u^{3}(t^{'}) \, dt^{'}
\end{array}
\label{Eqn:StrongSol}
\end{equation}

It is known \cite{lindsogge} that (\ref{Eqn:WaveEq}) is locally
well-posed for $s > \frac{1}{2}$ in $H^{s}(\mathbb{R}^{3}) \times
H^{s-1}(\mathbb{R}^{3})$ endowed with the standard norm $ \| (f,g)
\|_{H^{s} \times H^{s-1}} :=  \| f \|_{H^{s}} + \| g \|_{H^{s-1}}$.
Moreover the time of local existence does only depend on the norm of
the initial data $ \| (u_{0}, \, u_{1}) \|_{H^{s} \times H^{s-1}}$.

Now we turn our attention to the global well-posedness theory of
(\ref{Eqn:WaveEq}). In view of the above local well-posedness theory
and standard limiting arguments it suffices to establish an a priori
bound of the form

\begin{equation}
\begin{array}{ll}
\| \left( u(T), \partial_{t} u (T) \right) \|_{H^{s} \times H^{s-1}}
& \leq  C \left( s, \| u_{0} \|_{H^{s}}, \| u_{1} \|_{H^{s-1}}, T
\right)
\end{array}
\end{equation}
for all times $0 < T < \infty $ and all smooth-in-time
Schwartz-in-space solutions $(u,\partial_{t} u): [0, \, T] \times
\mathbb{R}^{3} \rightarrow \mathbb{R}$, where the right-hand side is
a finite quantity depending only on $s$, $ \| u_{0} \|_{H^{s}} $, $
\| u_{1} \|_{H^{s-1}}$ and $T$. Therefore in the sequel we shall
restrict attention to such smooth solutions.

The defocusing cubic wave equation (\ref{Eqn:WaveEq}) enjoys the
following energy conservation law

\begin{equation}
\begin{array}{ll}
E(u(t)) & := \frac{1}{2} \int_{\mathbb{R}^{3}} (\partial_{t}
u)^{2}(t,x) \, dx  + \frac{1}{2} \int_{\mathbb{R}^{3}} | D u
(t,x)|^{2} \, dx + \frac{1}{4} \int_{\mathbb{R}^{3}} u^{4}(t,x) \,
dx
\end{array}
\end{equation}
Combining this conservation law to the local well-posedness theory
we immediately have global well-posedness for (\ref{Eqn:WaveEq}) and
for $s=1$.

In this paper we are interested in studying global well-posedness of
(\ref{Eqn:WaveEq}) for data whose norm is below the energy norm, i.e
$s<1$. It is conjectured that (\ref{Eqn:WaveEq}) is globally
well-posed in $H^{s} (\mathbb{R}^{3}) \times H^{s-1}(\mathbb{R}^{3})
$ for $s > \frac{1}{2}$. The study of global existence for the
defocusing cubic wave equation has attracted the attention of many
researchers. Let us some mention some results for data
$(u_{0},u_{1})$ lying in a slightly different space than $H^{s}
\times H^{s-1}$ i.e $ \dot{H}^{s} \cap L^{4} \times \dot{H}^{s-1}$.
Here $\dot{H}^{s}$ is the standard homogeneous Sobolev space i.e the
completion of Schwartz functions $\mathcal{S} \left( \mathbb{R}^{3}
\right)$ with respect to the norm

\begin{equation}
\begin{array}{ll}
\| f \|_{\dot{H}^{s}} & := \| D^{s} f \|_{L^{2} \left(
\mathbb{R}^{3} \right)}
\end{array}
\end{equation}
Kenig, Ponce and Vega \cite{kenponcevega} were the first to prove
that (\ref{Eqn:WaveEq}) is globally well-posed for $ \frac{3}{4} <s
< 1$. They used the \emph{Fourier truncation method} discovered by
Bourgain \cite{bourg}. I. Gallagher and F. Planchon
\cite{gallagplanch}  proposed a different method to prove global
well-posedness for $1 > s > \frac{3}{4}$. H. Bahouri and Jean-Yves
Chemin \cite{bahchemin} proved global-wellposedness for
$s=\frac{3}{4}$ by using a non linear interpolation method and
logarithmic estimates from S. Klainermann and D. Tataru
\cite{klaintat}. Recently it was proved \cite{triroy} that the
defocusing cubic wave equation under spherically symmetric data is
globally well-posed in $H^{s} \times H^{s-1}$ for $1 > s >
\frac{7}{10} $. The main result of this paper is the following one

\begin{thm}
The defocusing cubic wave equation is globally well-posed in $ H^{s}
\times H^{s-1} $, $ 1 > s > \frac{13}{18} \simeq 0.722$. Moreover if
$s>\frac{13}{18}$ is close to $\frac{13}{18}$ then

\begin{equation}
\begin{array}{ll}
\| ( u(T), \, \partial_{t}u(T) ) \|^{2}_{H^{s} \times H^{s-1}} \leq
C \left( \| u_{0} \|_{H^{s}}, \| u_{1} \|_{H^{s-1}} \right)
T^{\frac{28s-18}{18s-13} +}
\end{array}
\label{Eqn:BoundHs}
\end{equation}
Here $C \left( \| u_{0} \|_{H^{s}}, \| u_{1} \|_{H^{s-1}}  \right)$
is a constant depending only on $\| u_{0} \|_{H^{s}}$ and $\| u_{1}
\|_{H^{s-1}}$. \label{Thm:Gwp}
\end{thm}
We set some notation that appear throughout the paper. Given $A,B$
positive number $A \lesssim B$ means that there exists a universal
constant $K$ such that $A \leq K B$. We say that $K_{0}$ is the
constant determined by the relation $A \lesssim B$ if $K_{0}$ is the
smallest $K$ such that $A \leq K B$ is true. We write $A \sim B$
when $A \lesssim B$ and $B \lesssim A$. $A << B$ denotes $ A \leq K
B$ for some universal constant $K < \frac{1}{100}$ . We also use the
notations $A+ = A + \epsilon$, $A++=A + 2 \epsilon$, $A- =A -
\epsilon$ and $A-- = A - 2 \epsilon$, etc. for some universal
constant $0 < \epsilon << 1$. We shall abuse the notation and write
$+$, $-$ for $0+$, $0-$ respectively. Let $\nabla$ denote the
gradient operator. If $J$ is an interval then $|J|$ is its size. Let
$I$ be the following multiplier

\begin{equation}
\begin{array}{ll}
\widehat{If}(\xi) & := m(\xi) \hat{f}(\xi)
\end{array}
\end{equation}
where $m(\xi): =  \eta \left( \frac{\xi}{N} \right)$, $\eta$ is a
smooth, radial, nonincreasing in $|\xi|$ such that

\begin{equation}
\begin{array}{ll}
\eta (\xi) & :=  \left\{
\begin{array}{l}
1, \, |\xi| \leq 1 \\
\left( \frac{1}{|\xi|} \right)^{1-s}, \, |\xi| \geq 2
\end{array}
\right.
\end{array}
\end{equation}
and $N>>1$ is a dyadic number playing the role of a parameter to be
chosen. We shall abuse the notation and write $m (|\xi|)$ for
$m(\xi)$, thus for instance $m(N)=1$. We denote by $E \left( I u(t)
\right)$ the so-called mollified energy

\begin{equation}
\begin{array}{ll}
E \left( Iu(t) \right) & := \frac{1}{2} \int_{\mathbb{R}^{3}} \left(
\partial_{t} I u(t,x) \right)^{2} \, dx + \frac{1}{2} \int_{\mathbb{R}^{3}} | D I u(t,x)|^{2} \, dx +
\frac{1}{4} \int_{\mathbb{R}^{3}} \left( I u(t,x) \right)^{4} \, dx
\end{array}
\label{Eqn:SmthNrj}
\end{equation}
The following result establishes the link between $ \| \left( u(T),
\partial_{t}u(T) \right) \|_{H^{s} \times H^{s-1}}$ and the
mollified energy $E(Iu)$ for a function $u$.

\begin{prop}{\cite{triroy}}
Let $T>0$. Then

\begin{equation}
\begin{array}{ll}
\| \left( u(T), \, \partial_{t} u(T) \right) \|^{2}_{ H^{s} \times
H^{s-1}} & \lesssim \| u_{0} \|^{2}_{H^{s}} + \left( T^{2}+1
\right)\sup_{t \in [0, \, T]} E \left( Iu(t) \right)
\end{array}
\label{Eqn:NrjEstim1}
\end{equation}
\label{Prop:NrjEst}
\end{prop}

We recall some basic results regarding the defocusing cubic wave
equation. Let $\lambda \in \mathbb{R}$ and $u_{\lambda}$ denote the
following function

\begin{equation}
\begin{array}{ll}
u_{\lambda}(t,x) & := \frac{1}{\lambda} u \left( \frac{t}{\lambda},
\frac{x}{\lambda} \right)
\end{array}
\end{equation}
If $u$ satisfies (\ref{Eqn:WaveEq}) with data $(u_{0},u_{1})$ then
$u_{\lambda}$ also satisfies (\ref{Eqn:WaveEq}) but with data $
\left( \frac{1}{\lambda} u_{0} \left( \frac{x}{\lambda} \right),
\frac{1}{\lambda^{2}} u_{1} \left( \frac{x}{\lambda} \right)
\right)$.

Now we recall the Strichartz estimates with derivative. These
estimates are proved in \cite{triroy} and follow from the standard
Strichartz estimates for the wave equation (\cite{ginebvelo},
\cite{lindsogge}).

\begin{prop} {\textbf{'' Strichartz estimates with derivative in $3$ dimensions''}}
Let $m \in [0, \, 1]$ and $ 0 \leq \tau < \infty$. If $u$ is a
strong solution to the IVP problem

\begin{equation}
\left\{
\begin{array}{l}
 \partial_{tt} u - \Delta   u = F \\
 u(0,x) = f(x) \\
 \partial_{t} u(0,x) =  g(x) \\
\end{array}
\right. \label{Eqn:LinearWave}
\end{equation}
then we have the $m$- Strichartz estimate with derivative

\begin{equation}
\begin{array}{l}
\| u \|_{L_{t}^{q}\left( [0, \, \tau] \right) L_{x}^{r}} + \|
\partial_{t} D^{-1} u \|_{L_{t}^{q} \left( [0, \, \tau] \right) L_{x}^{r}} + \| u \|_{L_{t}^{\infty} \left( [0, \, \tau] \right) \dot{H}^{m} } +
\| \partial_{t} u \|_{L_{t}^{\infty} \left( [0, \, \tau] \right)
\dot{H}^{m-1} }  \\
\lesssim \| f \|_{\dot{H}^{m}} + \| g \|_{\dot{H}^{m-1}} + \| F
\|_{L_{t}^{\tilde{q}} \left( [0, \, \tau] \right) L_{x}^{\tilde{r}}}
\end{array}
\label{Eqn:StrDer}
\end{equation}
under two assumptions

\begin{itemize}

\item $(q,r)$ lie in the set $\mathcal{W}$ of \textit{wave-admissible} points i.e

\begin{equation}
\begin{array}{ll}
\mathcal{W}  & := \left\{ (q, \, r): ( q, \, r ) \in (2, \,\infty ]
\times [2,\infty), \, \frac{1}{q}+\frac{1}{r} \leq \frac{1}{2}
\right\}
\end{array}
\label{Eqn:StrCondition1}
\end{equation}

\item $(\tilde{q}, \tilde{r})$ lie in the dual set $\widetilde{\mathcal{W}}$ of $\mathcal{W}$ i.e

\begin{equation}
\begin{array}{ll}
\widetilde {\mathcal{W}}  & :=  \left\{ (\tilde{q}, \tilde{r}):
\frac{1}{\tilde{q}} + \frac{1}{q}=1, \, \frac{1}{r}+
\frac{1}{\tilde{r}}=1, \, (q,r) \in \mathcal{W} \right\}
\end{array}
\end{equation}

\item $(q,r,\tilde{q}, \tilde{r})$ satisfy the \textit{dimensional analysis}
conditions

\begin{equation}
\begin{array}{l}
\frac{1}{q} + \frac{3}{r} =  \frac{3}{2} -m
\end{array}
\label{Eqn:StrCondition2}
\end{equation}
and

\begin{equation}
\begin{array}{l}
\frac{1}{\tilde{q}} + \frac{3}{\tilde{r}} -2 = \frac{3}{2} -m
\end{array}
\label{Eqn:StrCondition3}
\end{equation}

\end{itemize}
\label{Prop:StriCha}
\end{prop}

Some variables frequently appear in this paper. We define them now.

We say that $(q,r)$ is a $m$-wave admissible pair if $0 \leq  m \leq
1$ and $(q,r)$ satisfy the two following conditions

\begin{itemize}
\item $(q,r) \in \mathcal{W}$

\item $ \frac{1}{q} + \frac{3}{r} =  \frac{3}{2}-m $

\end{itemize}

Let $J=[a, \,b]$ be an interval included in $[0, \, \infty)$. Given
a function $u$ we define $Z_{m,s}(J, \, u)$

\begin{equation}
\begin{array}{ll}
Z_{m,s}(J, \, u) & :=  \sup_{q,r} \left( \| D^{1-m} Iu
\|_{L_{t}^{q}(J)L_{x}^{r}} + \| D^{-m} \partial_{t} I u
\|_{L_{t}^{q}(J)L_{x}^{r}} \right)
\end{array}
\end{equation}
where the $\sup$ is taken over $m$-wave admissible $(q,r)$  and let

\begin{equation}
\begin{array}{ll}
Z(J, \, u) & :=  \sup_{m \in [0, \, 1)} Z_{m,s}(J, \, u)
\end{array}
\end{equation}
If $u$ satisfies (\ref{Eqn:StrongSol}) then for $t \in J$

\begin{equation}
\begin{array}{ll}
u(t) & = \ulJ(t) + \unlJ(t)
\end{array}
\label{Eqn:DecompLinNl}
\end{equation}
with $\ulJ$ denoting the linear part of $u$ adapted to $J$ i.e

\begin{equation}
\begin{array}{ll}
\ulJ(t) & := \cos{\left((t-a)D \right)} \, u(a) + \frac {\sin{
\left( (t-a)D \right)}}{D} \,
\partial_{t}u(a)
\end{array}
\end{equation}
and $\unlJ$ denoting its corresponding nonlinear part i.e

\begin{equation}
\begin{array}{ll}
\unlJ(t) & := - \int_{a}^{t} \frac{\sin{ \left( (t-t^{'})D
\right)}}{D} u^{3}(t^{'}) \, dt^{'}
\end{array}
\end{equation}

\vspace{2mm}

Some estimates that we establish throughout the paper require a
Paley-Littlewood decomposition. We set it up now. Let $\phi(\xi)$ be
a real, radial, nonincreasing function that is equal to $1$ on the
unit ball $\left\{ \xi \in \mathbb{R}^{3}: \, |\xi| \leq 1 \right\}$
and that that is supported on $\left\{ \xi \in \mathbb{R}^{3}: \,
|\xi| \leq 2 \right\}$. Let $\psi$ denote the function

\begin{equation}
\begin{array}{ll}
\psi(\xi) & := \phi(\xi) - \phi(2 \xi)
\end{array}
\end{equation}
If $M \in 2^{\mathbb{Z}}$ is a dyadic number we define the
Paley-Littlewood operators in the Fourier domain by

\begin{equation}
\begin{array}{ll}
\widehat{P_{\leq M} f}(\xi) & := \phi \left( \frac{\xi}{M} \right)
\hat{f}(\xi) \\
\widehat{P_{M} f}(\xi) & := \psi \left( \frac{\xi}{M} \right)
\hat{f}(\xi) \\
\widehat{P_{> M} f}(\xi) & := \hat{f}(\xi) - \widehat{P_{\leq M}
f}(\xi)
\end{array}
\end{equation}
Since $\sum_{M \in  2^{\mathbb{Z}}} \psi \left( \frac{\xi}{M}
\right)=1$ we have

\begin{equation}
\begin{array}{ll}
f & = \sum_{M \in 2^{\mathbb{Z}}} P_{M} f
\end{array}
\end{equation}

We conclude this introduction by giving the main ideas of the proof
of Theorem \ref{Thm:Gwp} and explaining how the paper is organized.
We are interested in finding an a priori upper bound of $ \| \left(
u(T), \, \partial_{t} u(T) \right) \|_{H^{s} \times H^{s-1}} $.
Proposition \ref{Prop:NrjEst} shows that it suffices to estimate
$\sup_{t \in [0, \, T]} E \left( I u(t) \right)$. The variation of
the mollified energy is expected to be slow. Therefore our strategy
is to estimate the supremum of the mollified energy by applying the
fundamental theorem of calculus. This is the $I$- method originally
invented by J. Colliander, M. Keel, G. Staffilani, H. Takaoka and
T.Tao in \cite{almckstt} to prove global well-posedness for
semilinear Schr\"odinger equations and for rough data and designed
in \cite{triroy} for the defocusing cubic wave equation. We divide
the whole interval $[0,T]$ into same length intervals. On each of
these subintervals we estimate the variation of the mollified energy
by performing a Paley-Littlewood decomposition and, roughly
speaking, by dividing the pieces of the solution supported on high
frequencies into their linear part adapted to the subinterval and
their corresponding nonlinear part. We prove in Section
\ref{sec:LocalGbBd} that we can locally control some quantities
depending on the solution and globally control other quantities
depending on its linear part. Kenig, Ponce and Vega
\cite{kenponcevega} observed that the nonlinear part of $u$ is
smoother than the linear part on high frequencies. In the same
spirit we prove a local inequality in Section \ref{sec:NonRegPart}
that brings out this fact. The variation of the mollified energy
comprises three types of terms. Some of them are only made up of the
nonlinear part of the solution: they are estimated by using the gain
of regularity found in Section \ref{sec:NonRegPart} and they are
locally small but globally large. Some other are composed of the
linear part of the solution: they are estimated by using the global
estimates found in Section \ref{sec:LocalGbBd} and they are locally
larger but globally smaller. The other ones are mixed terms and are
estimated by using the results of Sections \ref{sec:LocalGbBd} and
\ref{sec:NonRegPart}: we expect a combination of both effects. We
estimate in Section \ref{sec:AlmCon} the variation of the smoothed
energy on each of these subintervals. Then we iterate to cover the
whole interval. The upper bound of the total variation depends on
the size of the subintervals. This one plays the role of a parameter
to be chosen. By minimizing the upper bound we find the optimal
value that yields the sharpest estimate. This process is explained
in Section \ref{sec:PfGw57}: Theorem \ref{Thm:Gwp} follows.

\begin{rem}
If we had used the original $I$ method \cite{almckstt} we would have
obtained a $O \left( \frac{1}{N^{1-}} \right)$ increase of the
smoothed energy on time intervals of size one and we would have
found global well-posedness for $s>\frac{3}{4}$ (see
\cite{triroythesis}) . In this paper we prove that we have the same
increase but on time intervals of size larger than one and this is
why we beat $\frac{3}{4}$ \footnote{More precisely the size is $\sim
N^{\frac{1}{4}}$: see (\ref{Eqn:EstNrjProof}).}.
\end{rem}

\vspace{5mm}

$\textbf{Acknowledgements}:$ The author would like to thank his
advisor Terence Tao for introducing him to this topic and is
indebted to him for many helpful conversations and encouragement
during the preparation of this paper.

%Z(tau) to define%

\section{Proof of global well-posedness in $H^{s} \times H^{s-1}$, $1 > s > \frac{13}{18} $}
\label{sec:PfGw57}

In this section we prove the global existence of (\ref{Eqn:WaveEq})
in  $H^{s} \times H^{s-1}$, $1 > s > \frac{13}{18} $. Our proof
relies on some intermediate results that we prove in the next
sections. More precisely we shall show the following

\begin{prop}{\textbf{''Local and Global Boundedness''}}

Let $J=[a, \, b]$ be an interval included in $[0, \, \infty]$.
Assume that $u$ satisfies (\ref{Eqn:WaveEq}) and that

\begin{equation}
\begin{array}{ll}
\sup_{t \in J} E \left( I u(t) \right) &  \leq 2
\end{array}
\label{Eqn:BdCond1}
\end{equation}
Then

\begin{equation}
\begin{array}{ll}
Z \left(J, \, \ulJ \right) & \lesssim 1
\end{array}
\label{Eqn:GlobalBdulin}
\end{equation}
Moreover if $(q,r)$ are $m$- wave admissible then

\begin{equation}
\begin{array}{ll}
\| D^{1-m} I u \|_{L_{t}^{q}(J) L_{x}^{r}} & \lesssim \left( \maxJ
\right)^{\frac{1}{q}}
\end{array}
\label{Eqn:LocalBdu}
\end{equation}
and

\begin{equation}
\begin{array}{ll}
\| D^{1-m} I \unlJ \|_{L_{t}^{q}(J) L_{x}^{r}} & \lesssim \left(
\maxJ \right)^{\frac{1}{q}}
\end{array}
\label{Eqn:LocalBdulin}
\end{equation}

\label{prop:LocalBdLinNlin}
\end{prop}

\begin{prop}{\textbf{'' Local Gain of Regularity of the Nonlinear Term''}}
Let $J=[a, \, b]$ be an interval included in $[0, \, \infty)$ and
$u$ such that (\ref{Eqn:WaveEq}) and (\ref{Eqn:BdCond1}) hold. Then

\begin{equation}
\begin{array}{ll}
\| \partial_{t} I u^{nl,J} \|_{L_{t}^{6}(J) L_{x}^{3}} +  \| D I
u^{nl,J} \|_{L_{t}^{6}(J) L_{x}^{3}} & \lesssim
(\max{(1,|J|)})^{\frac{2}{3}}
\end{array}
\label{Eqn:NonLinReg}
\end{equation}
\label{prop:NonLinReg}
\end{prop}

\begin{prop}{\textbf{''Almost Conservation Law ''}}
Let $J=[a, \, b]$ be an interval included in $[0,\infty)$ and $u$
such that (\ref{Eqn:WaveEq}) and (\ref{Eqn:BdCond1}) hold. Then

\begin{equation}
\begin{array}{ll}
\left| \sup_{t \in J} E ( Iu(t) ) - E(Iu(a)) \right| & \lesssim
\max{ \left(
\begin{array}{l}
\frac{\left( \maxJ \right)^{\frac{1}{2}} }{N^{1-}},\frac
{\maxJ^{\frac{5}{2}}}{N^{2-}}
\end{array}
\right) }
\end{array}
\label{Eqn:EstNrj}
\end{equation}
\label{prop:EstNrj}
\end{prop}
For the remainder of the section we show that Proposition
\ref{prop:EstNrj} implies Theorem \ref{Thm:Gwp}.

Let $T>0$ and $N=N(T) >>1 $ be a parameter to be chosen later. There
are three steps to prove Theorem \ref{Thm:Gwp}.

\begin{enumerate}

\item \textbf{Scaling}. It was proved in \cite{triroy} that there
exists $ C_{0}=C_{0} \left( \| u_{0} \|_{H^{s}}, \| u_{1}
\|_{H^{s-1}} \right)$ such that if $\lambda$ satisfies

\begin{equation}
\begin{array}{ll}
\lambda & =  C_{0} N^{\frac{{2(1-s)}}{2s-1}}
\end{array}
\label{Eqn:UpperBdLambda}
\end{equation}
then

\begin{equation}
\begin{array}{ll}
E \left( I u_{\lambda}(0) \right)  & \leq  \frac{1}{2}
\end{array}
\label{Eqn:InitTruncNrjEst}
\end{equation}

\item \textbf{Boundedness of the mollified energy}. Let $F_{T}$ denote the following set

\begin{equation}
\begin{array}{ll}
F_{T}= \left\{ T^{'} \in [0, \,T]: \sup_{t \in [0, \, \lambda
T^{'}]} E \left( I u_{\lambda}(t) \right) \leq 1  \right\}
\end{array}
\end{equation}
with $\lambda$ defined in (\ref{Eqn:UpperBdLambda}). We claim that
$F_{T}$ is the whole set $[0, \, T]$ for $N=N(T)
>> 1$ to be chosen later. Indeed

\begin{itemize}

\item $F_{T} \neq \emptyset$ since $0 \in F_{T}$ by
(\ref{Eqn:InitTruncNrjEst}).

\item $F_{T}$ is closed by the dominated convergence theorem.

\item $F_{T}$ is open. Let $\widetilde{T^{'}} \in F_{T}$. By
continuity there exists $\delta > 0$ such that for every $T^{'} \in
\left( \widetilde{T^{'}} - \delta, \, \widetilde{T^{'}} + \delta
\right) \cap [0, \, T]$ we have

\begin{equation}
\begin{array}{ll}
\sup_{t \in [0, \, \lambda T^{'}]} E \left( I u_{\lambda}(t) \right)
& \leq 2
\end{array}
\label{Eqn:HypInducNrj}
\end{equation}
Assume that $\lambda T^{'} \leq 1$. Then by
(\ref{Eqn:InitTruncNrjEst}), (\ref{Eqn:HypInducNrj}) and by
Proposition \ref{prop:EstNrj} and we have

\begin{equation}
\begin{array}{ll}
\left| \sup_{t \in [0, \lambda T^{'}]} E ( Iu_{\lambda}(t) ) -
\frac{1}{2} \right| & \lesssim \frac{1}{N^{1-}}
\end{array}
\label{Eqn:EstNrj0}
\end{equation}
Now assume that $\lambda T^{'} >1$. We divide the interval
$[0,\lambda T^{'}]$ into subintervals $(J_{i})_{i \in [1..l]}$ such
that $|J_{1}|=...=|J_{l-1}| =\epsilon$, $ \lambda T^{'}  \geq
\epsilon >1$ to be determined and $|J_{l}| \leq \epsilon$. By
(\ref{Eqn:InitTruncNrjEst}), (\ref{Eqn:HypInducNrj}) and Proposition
\ref{prop:EstNrj} we have

\begin{equation}
\begin{array}{ll}
\left| \sup_{t \in [0, \lambda T^{'}]} E ( Iu_{\lambda}(t) ) -
\frac{1}{2} \right| & \lesssim \frac{ \lambda T^{'}}{\epsilon}
\left( \frac{\epsilon^{\frac{1}{2}}}{N^{1-}} +
\frac{\epsilon^{\frac{5}{2}} }{N^{2-}} \right)
\end{array}
\label{Eqn:EstNrjProof}
\end{equation}
We are seeking to minimize the right-hand side of
(\ref{Eqn:EstNrjProof}) with respect to $\epsilon$. If $\lambda
T^{'} >> N^{\frac{1}{2}} $ then choosing $\epsilon \sim
N^{\frac{1}{2}}$ we have

\begin{equation}
\begin{array}{ll}
\left| \sup_{t \in [0, \lambda T^{'}]} E ( Iu_{\lambda}(t) ) -
\frac{1}{2} \right| \lesssim \frac{\lambda T}{N^{\frac{5}{4}-}}
\end{array}
\label{Eqn:EstNrj1}
\end{equation}
Now if $\lambda T^{'} \lesssim N^{\frac{1}{2}}$ then letting
$\epsilon=\lambda T^{'} $ we have

\begin{equation}
\begin{array}{ll}
\left| \sup_{t \in [0, \lambda T^{'}]} E ( Iu_{\lambda}(t) ) -
\frac{1}{2} \right| \lesssim  \frac{1}{N^{\frac{3}{4}-}}
\end{array}
\label{Eqn:EstNrj2}
\end{equation}
Let $C_{0}$, $C_{1}$ and $C_{2}$ be the constants determined by
$\lesssim$ in (\ref{Eqn:EstNrj0}), (\ref{Eqn:EstNrj1}) and
(\ref{Eqn:EstNrj2}) respectively and let
$C=\max{(C_{0},C_{1},C_{2})}$. Since $s
> \frac{13}{18}$ we can always choose for every $T>0$ a $N=N(T)>>1$
such that

\begin{equation}
\begin{array}{ll}
C \max{  \left( \frac{1}{N^{1-}}, \, \frac{\lambda
T}{N^{\frac{5}{4}-}}, \, \frac{1}{N^{\frac{3}{4}-}} \right) }  \leq
\frac{1}{2}
\end{array}
\label{Eqn:CdtionN}
\end{equation}
With this choice of $N=N(T) >>1$ we have $\sup_{t \in [0, \, \lambda
T^{'}]} E(I u_{\lambda}(t)) \leq 1$.
\end{itemize}
Hence $F_{T}=[0, \, T]$ with $N=N(T) >>1$ satisfying
(\ref{Eqn:CdtionN}).

\item \textbf{Conclusion}. Following the $I$- method described
in \cite{almckstt}

\begin{equation}
\begin{array}{ll}
\sup_{t \in [0, \, T]} E \left( I u(t) \right) & \lesssim \lambda
\sup_{t \in [0, \, \lambda T]} E (I u_{\lambda}(t)) \\
& \lesssim \lambda
\end{array}
\label{Eqn:ComparNrjScale}
\end{equation}
Combining (\ref{Eqn:ComparNrjScale}) and Proposition
\ref{Prop:NrjEst} we have global well-posedness in $H^{s} \times
H^{s-1}$, $1> s > \frac{13}{18}$. Now let $T$ be large and let
$s>\frac{13}{18}$ be close to $\frac{13}{18}$. Then let $N$ such
that

\begin{equation}
\begin{array}{ll}
\frac{0.9}{2} \leq  C \frac{ \lambda T }{N^{\frac{5}{4}-}} \leq
\frac{1}{2}
\end{array}
\label{Eqn:DetN}
\end{equation}
Notice that (\ref{Eqn:CdtionN}) is satisfied with this choice of
$N$. We plug (\ref{Eqn:DetN}) into (\ref{Eqn:ComparNrjScale}) and we
apply Proposition \ref{Prop:NrjEst} to get (\ref{Eqn:BoundHs}).

\end{enumerate}

\section{Proof of ''Local and Global Boundedness''}
\label{sec:LocalGbBd}

In this section we prove Proposition \ref{prop:LocalBdLinNlin}. In
what follows we also assume that $J=[0, \, \tau]$: the reader can
check after reading the proof that the other cases come down to this
one. We slightly modify an argument in \cite{triroy}. We multiply
the $m$-Strichartz estimate with derivative (\ref{Eqn:StrDer}) by
$D^{1-m} I$ and we have

\begin{equation}
\begin{array}{ll}
Z_{m,s}(\tau, \, \ulJ) & \lesssim \| D I u_{0} \|_{L^{2}} + \| I
u_{1}
\|_{L^{2}} \\
& \lesssim 1
\end{array}
\end{equation}
This proves (\ref{Eqn:GlobalBdulin}).

Now let us prove (\ref{Eqn:LocalBdu}) and (\ref{Eqn:LocalBdulin}).
Notice that it suffices by (\ref{Eqn:GlobalBdulin}) and the triangle
inequality to prove (\ref{Eqn:LocalBdu}). We divide $[0,\tau]$ into
subintervals $(J_{k})_{k \in [1,...l]}$ such that
$|J_{1}|=...=|J_{l-1}|=\tau_{0}$ and $ |J_{l}| \leq \tau_{0}$
$\tau_{0}>0$ constant to be determined. By concatenation it suffices
to establish that $Z(J_{k},u) \lesssim 1$, $k \in [1,...,l]$. We
will prove the claim for $k=1$. By iteration it is also true for
$k>1$. There are two steps

\begin{itemize}

\item  \textbf{First Step} We assume that $m \leq s$. We multiply the $m$-Strichartz estimate with derivative
(\ref{Eqn:StrDer}) by $D^{1-m}I$ and we get from the fractional
Leibnitz rule, the H\"older in time and the H\"older in space
inequalities

\begin{equation}
\begin{array}{ll}
Z_{m,s}(\tau_{0},u) & \lesssim  \| D I u_{0} \|_{L^{2}} + \| I u_{1}
\|_{L^{2}}  + \| D^{1-m} I (uuu) \|_{L_{t}^{1}([0, \, \tau_{0}])
L_{x}^{\frac{6}{5-2m}}} \\
& \lesssim  1+ \| D^{1-m} I u \|_{L_{t}^{\infty}\left( [0,
\,\tau_{0}] \right) L_{x}^{\frac{6}{3-2m}}} \| u \|_{L_{t}^{2}
\left( [0, \, \tau_{0}] \right)
L_{x}^{6}}^{2} \\
& \lesssim  1+ Z_{m,s} (\tau_{0},u) \left( \tau_{0}^{\frac{1}{3}} \|
P_{\leq N} u \|_{L_{t}^{6} ( [0, \,\tau_{0}] ) L_{x}^{6}} +
\tau_{0}^{s-\frac{1}{2}} \| P_{>N} u \|_{L_{t}^{\frac{1}{1-s}} ( [0,
\,
\tau_{0}] ) L_{x}^{6}}  \right)^{2} \\
& \lesssim  1 + Z_{m,s}(\tau_{0}, \, u) \left(
\tau_{0}^{\frac{1}{3}} \| Iu \|_{L_{t}^{6} ( [0, \, \tau_{0}])
L_{x}^{6}} + \tau_{0}^{s-\frac{1}{2}} \frac{ \| D^{1-s} Iu
\|_{L_{t}^{\frac{1}{1-s}} ( [0, \, \tau_{0}] )
L_{x}^{6}}}{N^{1-s}}  \right)^{2}  \\
& \lesssim  1 + Z_{m,s}(\tau_{0}, \, u) \left(\tau_{0}^{\frac{1}{3}}
\| Iu \|_{L_{t}^{6} ([0, \tau_{0}]) L_{x}^{6}} +
\tau_{0}^{s-\frac{1}{2}}
\frac{Z_{s,s}(\tau_{0}, \, u)}{N^{1-s}}  \right)^{2} \\
& \lesssim 1 + Z_{m,s}(\tau_{0}, \, u) \left(\tau_{0}^{\frac{1}{2}}
\left( \sup_{t \in J} E ( I u (t) ) \right)^{\frac{1}{6}} +
\tau_{0}^{s-\frac{1}{2}} \frac{Z_{s,s}(\tau_{0}, \, u)}{N^{1-s}}
\right)^{2}
\end{array}
\label{Eqn:ResZmsu}
\end{equation}

Assume that $m=s$. If $\tau_{0}>0$ is small enough then after
applying a continuity argument to (\ref{Eqn:ResZmsu}) we get from
(\ref{Eqn:BdCond1})

\begin{equation}
\begin{array}{ll}
Z_{s,s}(\tau_{0}, \, u) & \lesssim 1
\end{array}
\label{Eqn:Zmeqsu}
\end{equation}
Assume that $m<s$. Then by (\ref{Eqn:ResZmsu}) and
(\ref{Eqn:Zmeqsu})

\begin{equation}
\begin{array}{ll}
Z_{m,s}(\tau_{0}, \, u) & \lesssim  1
\end{array}
\label{Eqn:Zminfsu}
\end{equation}

\item \textbf{Second step} We assume that $m>s$. By
(\ref{Eqn:ResZmsu}), (\ref{Eqn:Zmeqsu}) and (\ref{Eqn:Zminfsu}) we
have

\begin{equation}
\begin{array}{ll}
\| D^{1-r} I (uuu) \|_{L_{t}^{1} [0,\tau_{0}]
L_{x}^{\frac{6}{5-2r}}} & \lesssim Z_{r,s}(\tau_{0}, \, u) \left(
\tau_{0}^{\frac{1}{2}} \left( \sup_{t \in [0, \, \tau_{0}]} E (
Iu(t)) \right)^{\frac{1}{6}} +
\frac{ \tau^{s-\frac{1}{2}} Z_{s,s}(\tau_{0}, \, u)}{N^{1-s}}  \right)^{2} \\
& \lesssim  1
\end{array}
\label{Eqn:BdFcTerm}
\end{equation}
for $r \leq s$. The inequality

\begin{equation}
\begin{array}{ll}
\| D^{1-m} I (uuu) \|_{L_{t}^{1} ([0, \, \tau_{0}])
L_{x}^{\frac{6}{5-2m}}} & \lesssim \| D^{1-r} I (uuu) \|_{L_{t}^{1}
([0, \, \tau_{0}]) L_{x}^{\frac{6}{5-2r}}}
\end{array}
\label{Eqn:SobIneqm}
\end{equation}
follows from the application of Sobolev  homogeneous embedding. We
multiply the $m$-Strichartz estimate with derivative
(\ref{Eqn:StrDer}) by $D^{1-m}I$ and get from (\ref{Eqn:BdFcTerm})
and (\ref{Eqn:SobIneqm})

\begin{equation}
\begin{array}{ll}
Z_{m,s}(\tau_{0}, \, u) & \lesssim \| D I u_{0} \|_{L^{2}} + \| I
u_{1} \|_{L^{2}} + \| D^{1-m} I (uuu) \|_{L_{t}^{1}
\left( [0, \, \tau_{0}] \right) L_{x}^{\frac{6}{5-2m}}} \\
& \lesssim 1
\end{array}
\end{equation}

\end{itemize}

\section{Proof of ''Local Gain of regularity of the nonlinear part''}
\label{sec:NonRegPart}

In this section we prove (\ref{Eqn:NonLinReg}). In what follows we
also assume that $J=[0, \, \tau]$: the reader can check after
reading the proof that the other cases come down to that one. We get
from Proposition \ref{Prop:StriCha} and Proposition
\ref{prop:LocalBdLinNlin}

\begin{equation}
\begin{array}{ll}
\| \partial_{t} I u^{nl,J} \|_{L_{t}^{6}([0, \, \tau]) L_{x}^{3}} +
\| D I u^{nl,J} \|_{L_{t}^{6}([0, \, \tau]) L_{x}^{3}} & \lesssim \|
D I(uuu) \|_{L_{t}^{\frac{3}{2}}([0, \tau]) L_{x}^{\frac{6}{5}}}
\end{array}
\label{Eqn:Gainunl1}
\end{equation}
But

\begin{equation}
\begin{array}{ll}
\| D I(uuu) \|_{L_{t}^{\frac{3}{2}}([0, \tau]) L_{x}^{\frac{6}{5}}}
& \lesssim \| D I u \|_{L_{t}^{\infty}([0, \tau]) L_{x}^{2}} \| u
\|^{2}_{L_{t}^{3} ([0, \tau]) L_{x}^{6}} \\
& \lesssim \| D I u \|_{L_{t}^{\infty}([0, \tau]) L_{x}^{2}} \left(
\| P_{< N} u \|_{L_{t}^{3} ([0, \tau]) L_{x}^{6}} +  \| P_{ \geq N}u
\|_{L_{t}^{3}([0, \tau]) L_{x}^{6}} \right)^{2} \\
& \lesssim  \| D I u \|_{L_{t}^{\infty}([0, \tau]) L_{x}^{2}} \left(
\tau^{\frac{1}{3}} \| D^{1-1} I u \|_{L_{t}^{\infty} ([0, \tau])
L_{x}^{6}} + \frac{ \| D^{1-\frac{2}{3}} I u \|_{L_{t}^{3} ([0,
\tau]) L_{x}^{6}}}{N^{\frac{1}{3}}}
\right)^{2} \\
& \lesssim \tau^{\frac{2}{3}}
\end{array}
\label{Eqn:Gainunl2}
\end{equation}
Combining (\ref{Eqn:Gainunl1}) and (\ref{Eqn:Gainunl2}) we get

\begin{equation}
\begin{array}{ll}
\| \partial_{t} I u^{nl,J} \|_{L_{t}^{6}([0, \, \tau]) L_{x}^{3}} +
\| D I u^{nl,J} \|_{L_{t}^{6}([0, \, \tau]) L_{x}^{3}} & \lesssim
\tau^{\frac{2}{3}}
\end{array}
\end{equation}

\section{Proof of ''Almost conservation law''}
\label{sec:AlmCon}

Let $J=[a, \, b]$ be an interval included in $[0,\infty)$ and $u$
such that (\ref{Eqn:WaveEq}) and (\ref{Eqn:BdCond1}) hold. Let $\tau
\in J$.  Then the Plancherel formula and the fundamental theorem of
calculus yield

\begin{equation}
\begin{array}{ll}
 \left| E ( Iu(\tau) ) - E ( Iu(a)) \right| & =
 \left| \int^{\tau}_{a} \int_{\xi_{1} + ... + \xi_{4} = 0}
\mu(\xi_{2},\xi_{3},\xi_{4}) \widehat{\partial_{t} I u}(t,\xi_{1})
\widehat{I u}(t,\xi_{2}) \widehat{I u} (t,\xi_{3}) \widehat{I u}
(t,\xi_{4}) \, d\xi_{2}...d \xi_{4} dt \right|
\end{array}
\label{Eqn:FstExpNrj}
\end{equation}
with

\begin{equation}
\begin{array}{ll}
\mu(\xi_{2},\xi_{3},\xi_{4}) & := 1 - \frac{m(\xi_{2}+ \xi_{3} +
\xi_{4})}{m(\xi_{2}) m(\xi_{3}) m(\xi_{4})}
\end{array}
\label{Eqn:DfnMu}
\end{equation}
We perform a Paley-Littlewood decomposition to estimate the right
hand side of (\ref{Eqn:FstExpNrj}). Let $u_{i}:=P_{N_{i}} u$ and let
$X$ denote the following number

\begin{equation}
\begin{array}{ll}
X & :=  \left| \int^{\tau}_{a} \int_{\xi_{1}+...+\xi_{4}=0}
\mu(\xi_{2},\xi_{3},\xi_{4}) \widehat{\partial_{t} I
u_{1}}(t,\xi_{1}) \widehat{I u_{2}} (t,\xi_{2}) \widehat{I u_{3}}
(t,\xi_{3}) \widehat{I u_{4}}(t,\xi_{4}) \, d\xi_{2}...d \xi_{4} dt
\right|
\end{array}
\label{Eqn:DfnXLocPlus}
\end{equation}
The strategy to estimate $X$ is explained in \cite{morckstt},
\cite{triroy}. We recall the main steps.

\emph{Overview of the strategy}.
\begin{enumerate}

\item \textbf{First step} We seek a pointwise bound of the symbol

\begin{equation}
\begin{array}{ll}
\left| \mu(\xi_{2},\xi_{3}, \xi_{4}) \right| & \leq
B(N_{2},N_{3},N_{4})
\end{array}
\end{equation}
Then  for some $A \subset \{1, \, ..., 4 \} $ to be chosen we
decompose for every $i \in A$ $u_{i}$ into its linear part
$u_{i}^{l}:=P_{N_{i}} u^{l,J}$ and its nonlinear part
$u_{i}^{nl}:=P_{N_{i}} u^{nl,J}$ and after expansion we need to
evaluate expressions of the form

\begin{equation}
\begin{array}{ll}
Y & := \left| \int^{\tau}_{a} \int_{\xi_{1}+...+\xi_{4}=0}
\mu(\xi_{2},\xi_{3},\xi_{4}) \widehat{\partial_{t} I
v_{1}}(t,\xi_{1}) \widehat{I v_{2}}(t,\xi_{2})... \widehat{I
v_{4}}(t,\xi_{4}) \, d\xi_{2}...d \xi_{4} dt \right|
\end{array}
\end{equation}
with $v_{j}$, $j \in \{1, \, ..., \,4 \}$ denoting $u_{j}^{nl}$ or
$u_{j}^{l}$ or $u_{j}$ \footnote{the value of  $v_{j}$ depends on
the choice of $A$}. We get from the Coifman-Meyer theorem
(\cite{coifmey}, p179)

\begin{equation}
\begin{array}{ll}
Y & \lesssim B(N_{2}, N_{3}, N_{4}) \|
\partial_{t} I v_{1} \|_{L_{t}^{p_{1}}(J) L_{x}^{q_{1}}}
\| I v_{2} \|_{L_{t}^{p_{2}}(J) L_{x}^{q_{2}}}...\| I v_{4}
\|_{L_{t}^{p_{4}}(J) L_{x}^{q_{4}}}
\end{array}
\label{Eqn:HoldXpm}
\end{equation}
with $(p_{j},q_{j})$, $j \in \{2, \, ..., \,4 \}$  such that $p_{j}
\in [1, \infty]$, $q_{j} \in (1, \infty)$, $\sum_{j=1}^{4}
\frac{1}{p_{j}}=1$, $\sum_{j=1}^{4} \frac{1}{q_{j}}=1$,
$(p_{j},q_{j})$ $m_{j}$-wave admissible for some $m_{j}^{'} \,s$
such that $0 \leq m_{j} < 1$ and $\frac{1}{p_{j}} +\frac{1}{q_{j}}
=\frac{1}{2}$ \footnote{In other words $(p_{j},q_{j}) = \left(
\frac{2}{m_{j}}, \, \frac{2}{1-m_{j}} \right)$}.

\item \textbf{Second Step} We use the following Bernstein inequalities

\begin{equation}
\begin{array}{ll}
\| I v_{j} \|_{L_{t}^{p_{j}}(J) L_{x}^{q_{j}}} & \lesssim
N_{j}^{m_{j}-1} \| D^{1-m_{j}} v_{j}
\|_{L_{t}^{p_{j}}(J) L_{x}^{q_{j}}} \\
\| \partial_{t} I v_{1} \|_{L_{t}^{p_{1}}(J) L_{x}^{q_{1}}} &
\lesssim  N_{1}^{m_{1}} \| D^{-m_{1}} \partial_{t} I
v_{1} \|_{L_{t}^{p_{1}}(J) L_{x}^{q_{1}}} \\
\| I v_{j} \|_{L_{t}^{6}(J) L_{x}^{3}} & \lesssim \frac{1}{N_{j}} \|
D I v_{j} \|_{L_{t}^{6}(J) L_{x}^{3}}
\end{array}
\label{Eqn:1}
\end{equation}
We plug (\ref{Eqn:1})  into (\ref{Eqn:HoldXpm}).

\item \textbf{Third step} The series must be summable. Therefore in some cases we might create
$N_{1}^{\pm}$, $N_{j}^{\pm}$ for some $j^{'}s$ by considering slight
variations  $(p_{1} \pm , \, q_{1} \pm)$, $(p_{j} \pm , \, q_{j}
\pm) \in [1, \, \infty] \times (1, \, \infty)$ of $(p_{1}, \,
q_{1})$, $(p_{j},q_{j})$ that are $m_{1} \, \pm$, $m_{j} \, \pm$ -
wave admissible and such that $\frac{1}{p_{1} \pm} + \frac{1}{q_{1}
\pm} =\frac{1}{2}$, $\frac{1}{p_{j} \pm} + \frac{1}{q_{j} \pm}
=\frac{1}{2}$ respectively. For instance if we create slight
variations $(p_{1}+, q_{1}+)$, $(p_{j}+,q_{j}+)$ of $(p_{1}+,
q_{1})$, $(p_{j},q_{j})$ respectively then we get from Bernstein and
H\"older in time inequalities

\begin{equation}
\begin{array}{ll}
\| I v_{j} \|_{L_{t}^{p_{j}+}(J) L_{j}^{q_{j}-}} & \lesssim
N_{j}^{-} N_{j}^{m_{j}-1} \| D^{1-(m_{j}-)} I v_{j}
\|_{L_{t}^{p_{j}+}(J) L_{x}^{q_{j}-}} \\
\| \partial_{t} I v_{1} \|_{L_{t}^{p_{1}+}(J) L_{x}^{q_{1}-}} &
\lesssim N_{1}^{-} N_{1}^{m_{1}} \| D^{-(m_{1}-)} \partial_{t} I
v_{1} \|_{L_{t}^{p_{1}+}(J) L_{x}^{q_{1}-}} \\
\| I v_{j} \|_{L_{t}^{6-}(J) L_{x}^{3+}} & \lesssim
\frac{N_{j}^{+}}{N_{j}} \| D I v_{j} \|_{L_{t}^{6}(J) L_{x}^{3}} \\
\| \partial_{t} I v_{1} \|_{L_{t}^{6-}(J) L_{x}^{3+}} & \lesssim
N_{1}^{+} \| \partial_{t}  I v_{1} \|_{L_{t}^{6}(J) L_{x}^{3}}
\end{array}
\end{equation}
It was proved \cite{triroy} that the following inequality holds
\footnote{More precisely  $\| I v_{j} \|_{L_{t}^{\frac{2}{1-
\epsilon}} L_{x}^{\frac{2}{\epsilon}}} \lesssim N_{j}^{\epsilon} \|
D^{1-(1- \epsilon^{'})} I v_{j} \|_{L_{t}^{\frac{2}{1-
\epsilon^{'}}} L_{x}^{\frac{2}{\epsilon^{'}}}}$ with $\epsilon^{'}=5
\epsilon$. }

\begin{equation}
\begin{array}{ll}
\| I v_{j} \|_{L_{t}^{2+}(J) L_{x}^{\infty-}} & \lesssim N_{j}^{+}
\| D^{1-(1-)} I v_{j} \|_{L_{t}^{2+}(J) L_{x}^{\infty-}}
\end{array}
\label{Eqn:IndirectNk}
\end{equation}
by using the localization in time to our advantage. The creation of
$N_{j}^{+}$ allows to make the summation with respect to $N_{j}$
whenever $N_{j} < 1$.

\end{enumerate}
This ends the overview of the strategy. \vspace{5 mm}

Let us get back to the proof. By symmetry we may assume that $N_{2}
\geq N_{3} \geq N_{4}$. Let $N_{1}^{*}$,..., $N_{4}^{*}$ be the four
numbers $N_{1}$,...,$N_{4}$ in order so that $N_{1}^{*} \geq
N_{2}^{*} \geq N_{3}^{*} \geq N_{4}^{*}$. We can assume that
$N_{1}^{*} \gtrsim N$ since if not the multiplier $\mu$ of $X$
vanishes and $X=0$. We can also assume that $N_{1}^{*} \sim
N_{2}^{*}$ since if not the convolution constraint
$\xi_{1}+...+\xi_{4}=0$ imposes $X=0$. There are three cases

\begin{itemize}

\item  \textbf{Case} $\mathbf{1}$: $N_{1}^{*}=N_{2}$ and $N_{2}^{*}=N_{1}$

We write  $u_{i}=u_{i}^{l,J} + u_{i}^{nl,J}$, $i \in \{1,2 \}$. We
need to estimate

\begin{equation}
\begin{array}{ll}
X_{1} & = \left| \int_{0}^{\tau} \int_{\xi_{1}+..+\xi_{4}=0}
\mu(\xi_{2},...,\xi_{4})
\widehat{\partial_{t}Iu_{1}^{l,J}}(t,\xi_{1})
\widehat{Iu_{2}^{l,J}}(t,\xi_{2}) \widehat{Iu_{3}}(t,\xi_{3})
\widehat{Iu_{4}}(t,\xi_{4}) \,  d\xi_{2}...d \xi_{4} dt \right|
\end{array}
\end{equation}

\begin{equation}
\begin{array}{ll}
X_{2} & = \left| \int_{0}^{\tau} \int_{\xi_{1}+..+\xi_{4}=0}
\mu(\xi_{2},...,\xi_{4})
\widehat{\partial_{t}Iu_{1}^{l,J}}(t,\xi_{1})
\widehat{Iu_{2}^{nl,J}}(t,\xi_{2}) \widehat{Iu_{3}}(t,\xi_{3})
\widehat{Iu_{4}}(t,\xi_{4}) \,  d\xi_{2}...d \xi_{4} dt \right|
\end{array}
\end{equation}

\begin{equation}
\begin{array}{ll}
X_{3} & = \left| \int_{0}^{\tau} \int_{\xi_{1}+..+\xi_{4}=0}
\mu(\xi_{2},...,\xi_{4})
\widehat{\partial_{t}Iu_{1}^{nl,J}}(t,\xi_{1})
\widehat{Iu_{2}^{l,J}}(t,\xi_{2}) \widehat{Iu_{3}}(t,\xi_{3})
\widehat{Iu_{4}}(t,\xi_{4}) \,  d\xi_{2}...d \xi_{4} dt \right|
\end{array}
\end{equation}
and

\begin{equation}
\begin{array}{ll}
X_{4} & = \left| \int_{0}^{\tau} \int_{\xi_{1}+..+\xi_{4}=0}
\mu(\xi_{2},...,\xi_{4}) \widehat{\partial_{t} I
u_{1}^{nl,J}}(t,\xi_{1}) \widehat{Iu_{2}^{nl,J}}(t,\xi_{2})
\widehat{Iu_{3}}(t,\xi_{3}) \widehat{Iu_{4}}(t,\xi_{4}) \,
d\xi_{2}...d \xi_{4} dt \right|
\end{array}
\end{equation}

There are two subcases

\begin{itemize}

\item \textbf{Case} $\mathbf{1.a}$: $N_{3} \gtrsim N$

We have

\begin{equation}
\begin{array}{ll}
 \left| \mu \right| & \lesssim \frac{m(N_{1})}{ m(N_{2}) m(N_{3})
 m(N_{4})} \\
 & \lesssim \frac{1}{m(N_{3}) m(N_{4})}
\end{array}
\label{Eqn:EstMultCase1a}
\end{equation}

By (\ref{Eqn:GlobalBdulin}), (\ref{Eqn:LocalBdu}),
(\ref{Eqn:LocalBdulin}), (\ref{Eqn:IndirectNk}) and
(\ref{Eqn:EstMultCase1a}) we have

\begin{equation}
\begin{array}{ll}
X_{1} & \lesssim \frac{1}{m(N_{3}) m(N_{4})} \| \partial_{t} I
u^{l,J}_{1} \|_{L_{t}^{2+}(J) L_{x}^{\infty -}} \| I u_{2}^{l,J}
\|_{L_{t}^{\infty}(J) L_{x}^{2}} \| I u_{3} \|_{L_{t}^{\infty -}(J)
L_{x}^{2 +}} \| I u_{4} \|_{L_{t}^{2+}(J) L_{x}^{\infty-}} \\
& \lesssim \frac{1}{m(N_{3}) m(N_{4})} N_{1}^{-} N_{1}
\frac{1}{N_{2}} \frac{N_{3}^{+}}{N_{3}}  N_{4}^{+} \| D^{-(1-)}
\partial_{t} I u^{l,J}_{1} \|_{L_{t}^{2+}(J) L_{x}^{\infty -}} \|
D^{1-0} I u_{2}^{l,J}
\|_{L_{t}^{\infty}(J) L_{x}^{2}} \\
& \| D^{1-(0+)} I u_{3} \|_{L_{t}^{\infty-}(J) L_{x}^{2+}} \| D^{1-(1-)} I u_{4} \|_{L_{t}^{2+}(J) L_{x}^{\infty-}} \\
& \lesssim \left(\maxJ \right)^{\frac{1}{2}} \frac{ N_{2}^{--}
N_{4}^{+}}{N^{1-}}
\end{array}
\end{equation}
Similarly

\begin{equation}
\begin{array}{ll}
X_{2} & \lesssim \frac{1}{m(N_{3}) m(N_{4})} \| \partial_{t} I
u^{l,J}_{1} \|_{L_{t}^{2+}(J) L_{x}^{\infty -}} \| I u_{2}^{nl,J}
\|_{L_{t}^{\infty}(J) L_{x}^{2}} \| I u_{3} \|_{L_{t}^{\infty-}(J)
L_{x}^{2+}} \| I u_{4} \|_{L_{t}^{2+}(J) L_{x}^{\infty-}} \\
& \lesssim \left(\maxJ \right)^{\frac{1}{2}} \frac{ N_{2}^{--}
N_{4}^{+}}{N^{1-}}
\end{array}
\end{equation}
We have

\begin{equation}
\begin{array}{ll}
X_{3} & \lesssim \frac{1}{m(N_{3}) m(N_{4})} \| \partial_{t} I
u^{nl,J}_{1} \|_{L_{t}^{\infty }(J) L_{x}^{2}} \| I u_{2}^{l,J}
\|_{L_{t}^{2+}(J) L_{x}^{\infty-}} \| I u_{3} \|_{L_{t}^{\infty-}(J)
L_{x}^{2+}} \| I u_{4} \|_{L_{t}^{2+}(J) L_{x}^{\infty-}} \\
& \lesssim \frac{1}{m(N_{3}) m(N_{4})} N_{2}^{-}
\frac{N_{3}^{+}}{N_{3}} N_{4}^{+} \| D^{-0}
\partial_{t} I u_{1}^{nl,J} \|_{L_{t}^{\infty}(J) L_{x}^{2}} \|
D^{1-(1-)} I u_{2}^{l,J} \|_{L_{t}^{2+}(J) L_{x}^{\infty-}} \\
& \| D^{1-(0+)} I u_{3} \|_{L_{t}^{\infty}(J) L_{x}^{2}} \|
D^{1-(1-)} I u_{4} \|_{L_{t}^{2+}(J) L_{x}^{\infty-}} \\
& \lesssim (\maxJ)^{\frac{1}{2}} \frac{ N_{2}^{--}
N_{4}^{+}}{N^{1-}}
\end{array}
\end{equation}
As for $X_{4}$ we make further decompositions. We write
$u_{3}=u_{3}^{l,J}+u_{3}^{nl,J}$ and we need to estimate

\begin{equation}
\begin{array}{ll}
X_{4,1} & = \left| \int_{a}^{\tau} \int_{\xi_{1}+..+\xi_{4}=0}
\mu(\xi_{2},...,\xi_{4}) \widehat{\partial_{t} I
u_{1}^{nl,J}}(t,\xi_{1}) \widehat{Iu_{2}^{nl,J}}(t,\xi_{2})
\widehat{Iu_{3}^{l,J}}(t,\xi_{3}) \widehat{Iu_{4}}(t,\xi_{4}) \,
d\xi_{2}...d \xi_{4} dt \right|
\end{array}
\end{equation}
and

\begin{equation}
\begin{array}{ll}
X_{4,2} & = \left| \int_{a}^{\tau} \int_{\xi_{1}+..+\xi_{4}=0}
\mu(\xi_{2},...,\xi_{4}) \widehat{\partial_{t} I
u_{1}^{nl,J}}(t,\xi_{1}) \widehat{Iu_{2}^{nl,J}}(t,\xi_{2})
\widehat{Iu_{3}^{nl,J}}(t,\xi_{3}) \widehat{Iu_{4}}(t,\xi_{4}) \,
d\xi_{2}...d \xi_{4} dt \right|
\end{array}
\end{equation}
We have

\begin{equation}
\begin{array}{ll}
X_{4,1} & \lesssim \frac{1}{m(N_{3}) m(N_{4})} \| \partial_{t} I
u_{1}^{nl,J} \|_{L_{t}^{\infty-}(J) L_{x}^{2+}} \| I u_{2}^{nl,J}
\|_{L_{t}^{\infty}(J) L_{x}^{2}} \| I u_{3}^{l,J} \|_{L_{t}^{2+}(J)
L_{x}^{\infty -}} \| I u_{4} \|_{L_{t}^{2+}(J) L_{x}^{\infty -}} \\
& \frac{1}{m(N_{3}) m(N_{4})}  N_{1}^{+} \frac{1}{N_{2}} N_{3}^{+}
N_{4}^{+} \| D^{-(0+)} \partial_{t} I u_{1}^{nl,J}
\|_{L_{t}^{\infty-}(J) L_{x}^{2+}} \| D I u_{2}^{nl,J}
\|_{L_{t}^{\infty}(J) L_{x}^{2}} \\
& \|D^{1-(1-)} I u_{3}^{l,J} \|_{L_{t}^{2+}(J) L_{x}^{\infty -}} \| D^{1-(1-)} I u_{4} \|_{L_{t}^{2+}(J) L_{x}^{\infty -}} \\
& \lesssim (\maxJ)^{\frac{1}{2}} \frac{ N_{2}^{--}
N_{4}^{+}}{N^{1-}}
\end{array}
\end{equation}
and by (\ref{Eqn:NonLinReg}) we have

\begin{equation}
\begin{array}{ll}
X_{4,2} & \lesssim \frac{1}{m(N_{3}) m(N_{4})} \| \partial_{t} I
u_{1}^{nl,J} \|_{L_{t}^{6-}(J) L_{x}^{3+}} \| I u_{2}^{nl,J}
\|_{L_{t}^{6}(J) L_{x}^{3}} \| I u_{3}^{nl,J} \|_{L_{t}^{6}(J)
L_{x}^{3}} \| I u_{4} \|_{L_{t}^{2+}(J) L_{x}^{\infty -}} \\
& \lesssim \frac{1}{m(N_{3}) m(N_{4})} N_{1}^{+} \frac{1}{N_{2}}
\frac{1}{N_{3}} N_{4}^{+} \| \partial_{t} I u_{1}^{nl,J}
\|_{L_{t}^{6}(J) L_{x}^{3}} \| D I u_{2}^{nl,J} \|_{L_{t}^{6}(J)
L_{x}^{3}} \\
& \|D I u_{3}^{nl,J} \|_{L_{t}^{6}(J) L_{x}^{3}} \|
D^{1-(1-)} I u_{4} \|_{L_{t}^{2+}(J) L_{x}^{\infty -}} \\
& \lesssim (\maxJ)^{\frac{5}{2}} \frac{N_{2}^{--} N_{4}^{+}}{N^{2-}}
\end{array}
\end{equation}

\item \textbf{Case} $\mathbf{1.b}$: $N_{3}<<N$
By the Mean Value Inequality

\begin{equation}
\begin{array}{ll}
|\mu| & \lesssim \frac{ |\nabla m(\xi_{2})| |\xi_{3} +
\xi_{4}|}{m(\xi_{2})} \\
& \lesssim \frac{N_{3}}{N_{2}}
\end{array}
\label{Eqn:EstMultCase1b}
\end{equation}
Now if for $X_{1}$, $X_{2}$, $X_{3}$, $X_{4,1}$ and $X_{4,2}$ we
apply the same procedure to that of Case $1.a$ and if use
(\ref{Eqn:EstMultCase1b}) we see that the factor $N_{3}^{\alpha}$
that appears always satisfies $\alpha \geq 0$ and consequently is
comparable to $N^{\alpha}$. Therefore the results are the same. For
instance

\begin{equation}
\begin{array}{ll}
X_{1} & \lesssim \frac{N_{3}}{N_{2}} \| \partial_{t} I u^{l,J}_{1}
\|_{L_{t}^{2+}(J) L_{x}^{\infty -}} \| I u_{2}^{l,J}
\|_{L_{t}^{\infty}(J) L_{x}^{2}} \| I u_{3} \|_{L_{t}^{\infty -}(J)
L_{x}^{2+}} \| I u_{4} \|_{L_{t}^{2+}(J) L_{x}^{\infty-}} \\
& \lesssim \frac{N_{3}}{N_{2}} N_{1}^{-} N_{1} \frac{1}{N_{2}}
\frac{N_{3}^{+}}{N_{3}}  N_{4}^{+} \| D^{-(1-)}
\partial_{t} I u^{l,J}_{1} \|_{L_{t}^{2+}(J) L_{x}^{\infty -}} \|
D^{1-0} I u_{2}^{l,J}
\|_{L_{t}^{\infty }(J) L_{x}^{2}} \\
& \| D^{1-(0+)} I u_{3} \|_{L_{t}^{\infty}(J) L_{x}^{2}} \| D^{1-(1-)} I u_{4} \|_{L_{t}^{2+}(J) L_{x}^{\infty-}} \\
& \lesssim \left(\maxJ \right)^{\frac{1}{2}} \frac{ N_{2}^{--}
N_{4}^{+}}{N^{1-}}
\end{array}
\end{equation}
and here the factor $ N_{3}^{+}\frac{N_{3}}{N_{3}}=N_{3}^{+}$
appears.

\end{itemize}

\item \textbf{Case} $\mathbf{2}$: $N_{1}^{*}=N_{1}$ and $N_{2}^{*}=N_{2}$

Since $N_{1} \sim N_{2}$ then this case boils down to the previous
one.

\item \textbf{Case} $\mathbf{3}$: $N_{1}^{*}=N_{2}$ and $N_{2}^{*}=N_{3}$

We write  $u_{i}=u_{i}^{l,J} + u_{i}^{nl,J}$, $i \in \{2,3 \}$. We
need to estimate

\begin{equation}
\begin{array}{ll}
X^{'}_{1} & = \left| \int_{a}^{\tau} \int_{\xi_{1}+..+\xi_{4}=0}
\mu(\xi_{2},...,\xi_{4}) \widehat{\partial_{t}Iu_{1}}(t,\xi_{1})
\widehat{Iu_{2}^{l,J}}(t,\xi_{2}) \widehat{Iu_{3}^{l,J}}(t,\xi_{3})
\widehat{Iu_{4}}(t,\xi_{4}) \, d\xi_{2}...d \xi_{4} dt \right|
\end{array}
\end{equation}

\begin{equation}
\begin{array}{ll}
X^{'}_{2} & = \left| \int_{a}^{\tau} \int_{\xi_{1}+..+\xi_{4}=0}
\mu(\xi_{2},...,\xi_{4}) \widehat{\partial_{t}Iu_{1}}(t,\xi_{1})
\widehat{Iu_{2}^{l,J}}(t,\xi_{2}) \widehat{Iu_{3}^{nl,J}}(t,\xi_{3})
\widehat{Iu_{4}}(t,\xi_{4}) \,  d\xi_{2}...d \xi_{4} dt \right|
\end{array}
\end{equation}

\begin{equation}
\begin{array}{ll}
X^{'}_{3} & = \left| \int_{a}^{\tau} \int_{\xi_{1}+..+\xi_{4}=0}
\mu(\xi_{2},...,\xi_{4}) \widehat{\partial_{t}Iu_{1}} (t,\xi_{1})
\widehat{Iu_{2}^{nl,J}}(t,\xi_{2}) \widehat{Iu_{3}^{l,J}}(t,\xi_{3})
\widehat{Iu_{4}}(t,\xi_{4}) \,  d\xi_{2}...d \xi_{4} dt \right|
\end{array}
\end{equation}
and

\begin{equation}
\begin{array}{ll}
X^{'}_{4} & = \left| \int_{a}^{\tau} \int_{\xi_{1}+..+\xi_{4}=0}
\mu(\xi_{2},...,\xi_{4}) \widehat{\partial_{t} I u_{1}}(t,\xi_{1})
\widehat{Iu_{2}^{nl,J}}(t,\xi_{2})
\widehat{Iu^{nl,J}_{3}}(t,\xi_{3}) \widehat{Iu_{4}}(t,\xi_{4}) \,
d\xi_{2}...d \xi_{4} dt \right|
\end{array}
\end{equation}
We have

\begin{equation}
\begin{array}{ll}
|\mu| & \lesssim \frac{m(N_{1})}{m(N_{2})m(N_{3}) m(N_{4})}
\end{array}
\label{Eqn:EstMultCase3}
\end{equation}
By (\ref{Eqn:GlobalBdulin}), (\ref{Eqn:LocalBdu}),
(\ref{Eqn:LocalBdulin}), (\ref{Eqn:IndirectNk}) and
(\ref{Eqn:EstMultCase3}) we have

\begin{equation}
\begin{array}{ll}
X_{1}^{'} & \lesssim \frac{m(N_{1})}{m(N_{2})m(N_{3}) m(N_{4})} \|
\partial_{t} I u_{1} \|_{L_{t}^{\infty-}(J) L_{x}^{2+}} \| I u_{2}^{l,J}
\|_{L_{t}^{\infty}(J) L_{x}^{2} } \| I u_{3}^{l,J} \|_{L_{t}^{2+}(J)
L_{x}^{\infty -}} \| I u_{4} \|_{L_{t}^{2+}(J) L_{x}^{\infty -}} \\
& \lesssim \frac{m(N_{1})}{m(N_{2})m(N_{3}) m(N_{4})} N_{1}^{+}
\frac{1}{N_{2}} N_{3}^{+} N_{4}^{+} \| D^{-(0+)} \partial_{t} I
u_{1} \|_{L_{t}^{\infty -}(J) L_{x}^{2+}} \| D I u_{2}^{l,J}
\|_{L_{t}^{\infty}(J) L_{x}^{2}} \\
 &  \| D^{1-(1-)} I u_{3}^{l,J}
\|_{L_{t}^{2+}(J) L_{x}^{\infty -}} \| D^{1-(1-)} I u_{4}
\|_{L_{t}^{2+}(J) L_{x}^{\infty -}} \\
& \lesssim (\maxJ)^{\frac{1}{2}} \frac{N_{1}^{+} N_{2}^{---}
N_{4}^{+} }{N^{1-}}
\end{array}
\end{equation}
and

\begin{equation}
\begin{array}{ll}
X_{2}^{'} & \lesssim \frac{m(N_{1})}{m(N_{2})m(N_{3}) m(N_{4})} \|
\partial_{t} I u_{1} \|_{L_{t}^{\infty-}(J) L_{x}^{2}} \| I u_{2}^{l,J}
\|_{L_{t}^{2+}(J) L_{x}^{\infty -} } \| I u_{3}^{nl,J}
\|_{L_{t}^{\infty}(J)
L_{x}^{2}} \| I u_{4} \|_{L_{t}^{2+}(J) L_{x}^{\infty -}} \\
& \lesssim \frac{m(N_{1})}{m(N_{2})m(N_{3}) m(N_{4})} N_{1}^{+}
N_{2}^{+} \frac{1}{N_{3}}  N_{4}^{+} \| D^{-(0+)} \partial_{t} I
u_{1} \|_{L_{t}^{\infty -}(J) L_{x}^{2+}} \| D^{1-(1-)} I
u_{2}^{l,J} \|_{L_{t}^{2+}(J) L_{x}^{\infty -}} \\
& \| D I u_{3}^{nl,J}
\|_{L_{t}^{\infty}(J) L_{x}^{2}} \| D^{1-(1-)} I u_{4} \|_{L_{t}^{2+}(J) L_{x}^{\infty -}} \\
& \lesssim (\maxJ)^{\frac{1}{2}} \frac{N_{1}^{+} N_{2}^{---}
N_{4}^{+} }{N^{1-}}
\end{array}
\end{equation}
Similarly since $N_{2} \sim N_{3}$ we have

\begin{equation}
\begin{array}{ll}
X_{3}^{'} & \lesssim \frac{m(N_{1})}{m(N_{2})m(N_{3}) m(N_{4})} \|
\partial_{t} I u_{1} \|_{L_{t}^{\infty-}(J) L_{x}^{2}} \| I u_{2}^{nl,J}
\|_{L_{t}^{\infty}(J) L_{x}^{2} }  \| I u_{3}^{l,J}
\|_{L_{t}^{2+}(J) L_{x}^{\infty -}} \| I u_{4}
\|_{L_{t}^{2+}(J) L_{x}^{\infty -}} \\
& \lesssim  (\maxJ)^{\frac{1}{2}} \frac{N_{1}^{+} N_{2}^{---}
N_{4}^{+} }{N^{1-}}
\end{array}
\end{equation}
As for $X^{'}_{4}$ we make further decompositions. We write
$u_{1}=u_{1}^{l,J}+u_{1}^{nl,J}$ and we need to estimate

\begin{equation}
\begin{array}{ll}
X^{'}_{4,1} & = \left| \int_{a}^{\tau} \int_{\xi_{1}+..+\xi_{4}=0}
\mu(\xi_{2},...,\xi_{4}) \widehat{\partial_{t} I
u^{l,J}_{1}}(t,\xi_{1}) \widehat{Iu_{2}^{nl,J}}(t,\xi_{2})
\widehat{Iu_{3}^{nl}}(t,\xi_{3}) \widehat{Iu_{4}}(t,\xi_{4}) \,
d\xi_{2}...d \xi_{4} dt \right|
\end{array}
\end{equation}
and

\begin{equation}
\begin{array}{ll}
X^{'}_{4,2} & = \left| \int_{a}^{\tau} \int_{\xi_{1}+..+\xi_{4}=0}
\mu(\xi_{2},...,\xi_{4}) \widehat{\partial_{t} I
u_{1}^{nl,J}}(t,\xi_{1}) \widehat{Iu_{2}^{nl,J}}(t,\xi_{2})
\widehat{Iu_{3}^{nl,J}}(t,\xi_{3}) \widehat{Iu_{4}}(t,\xi_{4}) \,
d\xi_{2}...d \xi_{4} dt \right|
\end{array}
\end{equation}
We have

\begin{equation}
\begin{array}{ll}
X^{'}_{4,1} & \lesssim \frac{m(N_{1})}{m(N_{2}) m(N_{3}) m(N_{4})}
\| \partial_{t} I u_{1}^{l,J} \|_{L_{t}^{2+}(J) L_{x}^{\infty-}} \|
I u_{2}^{nl,J} \|_{L_{t}^{\infty}(J) L_{x}^{2}} \| I u_{3}^{nl,J}
\|_{L_{t}^{\infty-}(J) L_{x}^{2+}} \| I u_{4} \|_{L_{t}^{2+}(J)
L_{x}^{\infty-}} \\
& \lesssim \frac{m(N_{1})}{m(N_{2}) m(N_{3}) m(N_{4})} N_{1}^{-}
N_{1} \frac{1}{N_{2}} \frac{N_{3}^{+}}{N_{3}} N_{4}^{+} \| D^{-(1-)}
\partial_{t} I u_{1}^{l,J} \|_{L_{t}^{2+}(J) L_{x}^{\infty-}} \| D I
u_{2}^{nl,J} \|_{L_{t}^{\infty}(J) L_{x}^{2}} \\
& \| D^{1-(0+)} I u_{3}^{nl,J} \|_{L_{t}^{\infty-}(J) L_{x}^{2+}} \|
D^{1-(1-)} I
u_{4} \|_{L_{t}^{2+}(J) L_{x}^{\infty-}} \\
& \lesssim (\maxJ)^{\frac{1}{2}} \frac{N_{1}^{+} N_{2}^{---}
N_{4}^{+} }{N^{1-}}
\end{array}
\end{equation}
and by (\ref{Eqn:NonLinReg}) we have

\begin{equation}
\begin{array}{ll}
X^{'}_{4,2} & \lesssim \frac{m(N_{1})}{m(N_{2}) m(N_{3}) m(N_{4})}
\|
\partial_{t} I u_{1}^{nl,J} \|_{L_{t}^{6 -}(J) L_{x}^{3+}} \| I
u_{2}^{nl,J} \|_{L_{t}^{6}(J) L_{x}^{3}} \| I u_{3}^{nl,J}
\|_{L_{t}^{6}(J) L_{x}^{3}} \| I u_{4} \|_{L_{t}^{2+}(J)
L_{x}^{\infty-}} \\
& \lesssim \frac{m(N_{1})}{m(N_{2}) m(N_{3}) m(N_{4})} N_{1}^{+}
\frac{1}{N_{2}} \frac{1}{N_{3}} N_{4}^{+} \|
\partial_{t} I u_{1}^{nl,J} \|_{L_{t}^{6}(J) L_{x}^{3}} \| D I
u_{2}^{nl,J} \|_{L_{t}^{6}(J) L_{x}^{3}} \\
& \| D I u_{3}^{nl,J} \|_{L_{t}^{6}(J) L_{x}^{3}} \| D^{1-(1-)} I
u_{4} \|_{L_{t}^{2+}(J) L_{x}^{\infty-}} \\
& \lesssim (\maxJ)^{\frac{5}{2}} \frac{N_{1}^{+} N_{2}^{---}
N_{4}^{+} }{N^{2-}}
\end{array}
\end{equation}

\end{itemize}

\end{document}